\documentclass[12pt]{amsart}
\usepackage{amssymb}
\usepackage{amsxtra}
\usepackage[mathscr]{eucal}
\usepackage{graphics}
\usepackage{epsfig}

\newtheorem{theorem}{Theorem}[section]
\newtheorem{lemma}[theorem]{Lemma}
\newtheorem{corollary}[theorem]{Corollary}

\theoremstyle{definition}

\newtheorem{example}[theorem]{Example}

\def\C{\mathbb{C}}
\def\Z{\mathbb{Z}}
\def\Q{\mathbb{Q}}
\def\K{\mathbb{K}}
\def\E{\mathbf{E}}
\def\H{\mathbf{H}}
\def\cI{\mathcal{I}}

\def\Nm{{\mathrm{Nm}}}

\begin{document}
\title{Prime Divisors of Sequences Associated to Elliptic Curves}
\author{Graham Everest and Igor E. Shparlinski}
\address{(GE) School of Mathematics, University of East Anglia,
Norwich NR4 7TJ, UK}
\email{g.everest@uea.ac.uk}
\address{(IS) Department of Computing, Macquarie University,
NSW 2109, Australia}
\email{igor@comp.mq.edu.au}

\begin{abstract}
We consider the primes which divide the denominator of the
x-coordinate of a sequence of rational points on an elliptic
curve. It is expected that for every sufficiently large value of
the index, each term should be divisible by a primitive prime
divisor, one that has not appeared in any earlier term. Proofs of
this are known in only a few cases. Weaker results in the general
direction are given, using a strong form of Siegel's Theorem and
some congruence arguments. Our main result is applied to the study
of prime divisors of Somos sequences.
\end{abstract}
\maketitle

\section{Introduction}

Let $\E$ denote the elliptic
curve in generalized Weierstrass form,
\begin{equation}\label{model}
y^2+a_1xy+a_3y=x^3+a_2x^2+a_4x+a_6,
\end{equation}
defined over $\Q$, the rational numbers. We assume that the
coefficients are chosen to lie in $\Z$ and equation (\ref{model})
defines a minimal model. For background, definitions and all
properties of elliptic curves used in this paper,
consult~\cite{SIL1} and~\cite{SIL2}. Let $\K$ denote an algebraic
number field of degree $d=[\K:\Q]$ over $\Q$. Throughout the
paper, $\E(\K)$ denotes the group of $\K$-rational points of $\E$
and $O$  denotes the point at infinity, the identity for the group
of $\K$-rational points. Suppose $P$ and $Q$ denote $\K$-rational
points, $P,Q\in \E(\K)$, with $P$ non-torsion. Assume that $nP + Q
\ne O$ for all positive integers $n$ and  write $nP+Q=(x_n,y_n)$.
The assumptions on $\E$ allow the factorization
\begin{equation}\label{definegenbeta}
(x_n)=(x(nP+Q))=\alpha_n/\beta_n,
\end{equation}
of the principal ideal $(x(nP+Q))$
into relatively prime integral ideals $\alpha_n$ and $\beta_n$.

In the rational case, we may take $\beta_n$ to be a positive
integer. Silverman~\cite{SIL3} proved that when $P$ is a rational
point and $Q=O$ the point at infinity, for all sufficiently large
$n$, $\beta_n$ has a {\it primitive divisor\/}, that is,  a
divisor of $\beta_n$ which is coprime to $\beta_m$ for all
positive integer $m<n$. In general, the expression {\it primitive
ideal divisor\/} of a term $\beta_n$ is used to describe an ideal
$\cI$ which divdes $\beta_n$ but no $\beta_m$ with $m<n$. Cheon
and Hahn~\cite{ch} have extended Silverman's result~\cite{SIL3} to
the algebraic case, showing that for all sufficiently large $n$,
$\beta_n$ has a primitive ideal divisor. The same conclusions are
drawn in~\cite{ek} under the assumption that $Q$ is an arbitrary
torsion point. Under the assumption that $Q=O$ and $P$ is the
image of a $\K$-rational point under a non-trivial isogeny, it is
shown in~\cite{ek} that for all sufficiently large $n$, $\beta_n$
has a composite primitive ideal divisor. Pheidas~\cite{phi} has
raised the exciting possibility that results of this kind could be
related to Hilbert's Tenth Problem over the rationals, which has
still not been settled.

Results about primitive divisors have a long and fine tradition for
certain sequences which satisfy a linear recurrence
relation, see~\cite{epsw}.
Consider the sequence $u_n=a^n-b^n$, where $a>b$ are
positive coprime integers. Zsigmondy's
theorem~\cite{Z} says that
$u_n$ always has a primitive divisor unless (i) $a=2, b=1$
and $n=6$ or (ii) $a+b=2^k$ and $n=2$.
In~\cite{sch}, Schinzel extended this result in
a number of directions. He
established a large class of indices for which
$u_n$ has a composite primitive divisor.
As an example of his results, when $u$ is the Mersenne sequence
($a=2$ and $b=1$), he proved $u_{4k}$ has
a composite primitive divisor for all odd $k>5$.
He also proved an algebraic version for primitive
ideal divisors of sequences of the form $u_n=\alpha^n-\beta^n$,
where $\alpha$ and $\beta$ are algebraic numbers with
$\alpha/\beta$ is not a root of unity, see~\cite{Sch} and~\cite{Sch2}. If
the
sequence is changed slightly to
$u_n=\gamma\alpha^n-\beta^n$, then a Zsigmondy-type
theorem is proved in~\cite{Sch} under the assumption
that $\gamma$ is a root of unity.
The Zsigmondy theorem extends to all Lucas and Lehmer sequences
and has recently been cast in a
very strong uniform manner.
Bilu, Hanrot and Voutier prove
in~\cite{bhv} that for any $n>30$, the $n$-th term of any Lucas
or Lehmer sequence has a primitive divisor.
For an integral ideal $\cI$ of $\K$,
let $\omega_\K(\cI)$ denote the number of prime
ideal divisors of $\cI$. Assuming $P$ is non-torsion guarantees
that all of the terms in
the sequence $\beta$
are non-zero.
Obviously, for any sequence of integral ideals $v$ which
satisfies a Zsigmondy-type theorem, the following lower
bound holds:
\begin{equation}\label{lower}
\omega_\K\left(\prod_{n=M+1}^{M+N}v_n\right) \geq N,
\end{equation}
for all sufficiently large $M$. In the absence of a full Zsigmondy
theorem for elliptic sequences, we still prove a similar lower bound.

Recall that $A \ll B$ and $B \gg A$ are equivalent to $A= O(B)$.
Throughout the paper, the implied constants  may dependent on $\E$ and the
points $P$ and $Q$.

\begin{theorem}\label{weakZ}Let $P$ and $Q$ denote algebraic
points on an elliptic curve, with $P$ non-torsion. With $\beta_n$ defined
by~(\ref{definegenbeta}),
for all sufficiently large $M$,
\begin{equation}\label{weakZlb}
\omega_\K\left(\prod_{n=M+1}^{M+N}\beta_n\right) \gg N.
\end{equation}
\end{theorem}

Very few results about the existence of primitive
(rational) divisors are known for the Lehmer-Pierce sequence
$\Nm(\alpha^n-\beta^n)$ (see~\cite{L,P}),
where $\Nm(z)$ denotes the usual field norm of $z \in \K$. In~\cite{eg},
some easy counter-examples are given but the general problem
remains very much open.
In the elliptic case, the same problem is studied in~\cite{ek},
with greater success, although a fully general understanding
lies some way off. Write $b_n$ for the ideal norm
\begin{equation}\label{defineb}
b_n=\Nm(\beta_n).
\end{equation}
An immediate corollary of Theorem~\ref{weakZ} follows. When
$\K=\Q$, write $\omega_\K=\omega_{\Q}=\omega$ in the usual
way.

\begin{corollary}
\label{cor1}
 With $b_n$ defined
by~(\ref{defineb}),
for all sufficiently large $M$,
$$
\omega\left(\prod_{n=M+1}^{M+N}b_n\right) \gg N.
$$
\end{corollary}

The proof of the Corollary is immediate because $\omega_\K
(\beta)\leq d\omega(\Nm (\beta))$ for any ideal $\beta$ from $\K$.

Finally, we include some discussion to indicate what the
true picture might be for divisors of sequences
$\beta$ as in~(\ref{definegenbeta}). Calculations support
our belief that for any given $w$, for all sufficiently
large $n$, $\beta_n$  has a primitive ideal divisor with
at least $w$ distinct prime ideal factors. For the moment,
the best result we can obtain in that direction is the
following.

\begin{theorem}\label{better} With the hypotheses
of Theorem~\ref{weakZ}, suppose $Q = O$.
Given any
$w>0$, there is an integer $q\ge 1$
such that
$$
  \omega_\K \left( \prod_{n=M+1}^{M+N}\beta_{qn}\right)
\geq wN.
$$
\end{theorem}

In the next two sections, we gather some tools then prove
Theorem~\ref{weakZ}. In the following section,
the proof
of Theorem~\ref{better} is given. Finally,we give an
application of Theorem~\ref{weakZ} to the study
of divisors of Somos
sequences.

Our thanks go to Helen King for helping to clarify the proof of
Theorem \ref{weakZ}.

\section{Congruences and the Growth of $b_n$}

The first lemma is a very strong form of Siegel's Theorem.

\begin{lemma}\label{growth}
Suppose $P$ denotes
a non-torsion $\K$-rational point of $\E$.
For all
sufficiently large $n$
$$
\log b_n = dhn^2 + O(n),
$$
where the constant
$h=\widehat h(P)$ is the global canonical height of the underlying point
$P$.
\end{lemma}

\begin{proof}[Proof of Lemma~\ref{growth}.]
Writing $H(R)$ for the naive height of an algebraic
point, then
$$H(R)=\prod_v\max \{1,|x(R)|_v\}^{\frac{1}{d}},
$$
where the product is taken  over all valuations
of $\K$, including
the archimedean or infinite valuations, which correspond to the
embeddings of $\K$ into $\C$. In our set-up, we can interpret this as
$$\log H(nP+Q)=\frac{1}{d}\sum _{v|\infty} \log \max \{1,|x_n|_v\}
+\frac{1}{d}\log b_n.
$$
The theory of heights gives an estimate for
$$
\log H(nP+Q)=\widehat h(nP+Q)+O(1),
$$
where $\widehat h(nP+Q)$ denotes the canonical height of $nP+Q$. Since
this is a positive definite quadratic form, the
bound in~(\ref{b1}) follows:
\begin{equation}\label{b1}
\sum _{v|\infty} \log \max \{1,|x_n|_v\}+ \log b_n=dhn^2+O(n).
\end{equation}

The estimate in Lemma~\ref{growth} follows from
an upper bound for $|x_n|_v$ for each
archimedean valuation $v$. When this quantity
is large it means $nP+Q$ is close to the point at
infinity in that valuation. On the complex torus, this means the
elliptic logarithm is close to zero. Thus bounds
from elliptic transcendence theory are applicable. We use
Th\'eor\`eme 2.1 in~\cite{dav} but see also~\cite{ST}
where an explicit version of David's Theorem appears on
page 20.
The nature of the bound is
\begin{equation}\label{b2}
\log |x_n|_v \ll \log n \log \log n,
\end{equation}
where the implied constant depends upon $v$, $\E$ and the points
$P$ and $Q$. Combining the estimates in~(\ref{b1}) and~(\ref{b2})
gives Lemma~\ref{growth}.
\end{proof}

For  integers  $M \ge 0$,  $N \ge 1$, and $\cI$ an integral
ideal,
denote by $T(M,N,\cI)$
the number of solutions of the congruence
$$
\beta_n  \equiv 0 \pmod \cI, \qquad M+1 \le n \le M+N.
$$

\begin{lemma}\label{propofT}
\label{le:Cong} If $P$ is not a torsion point
and $\cI$ is an integral ideal then for any integers
$M \ge 0$,  $N \ge 1$,
$$
T(M,N,\cI)  \ll \frac{N}{(\log \Nm(\cI))^{1/2}} + 1.
$$
\end{lemma}

\begin{proof}
In this proof we are going to use the arithmetic
of the curve $\E$ modulo $\cI$ so some care is needed.
It is known, see~\cite{SIL1}, that for any ideal $\wp^r$,
where $\wp$ denotes a prime ideal of $\K$,
the set of rational points with
denominator divisible by $\wp^r$, together with $O$,
forms a group.
The set of points with denominators
divisible by $\cI$ is the intersection of those groups for
all powers of prime ideals dividing $\cI$; hence it too is a group
and is closed under subtraction between two points.

Suppose there are two integers $n$ and $k$ with $1
\le n < n+k  \le N$
and $k \le N/(T(M,N,\cI)  -1)$ for which
$$
\beta_n  \equiv \beta_{n+k}  \equiv 0 \pmod \cI.
$$
By the opening remark, subtracting gives again a point
with denominator divisible by $\cI$, $kP\equiv O \pmod \cI$.
Because $P$ is
not a torsion point we conclude that $b_k \ne 0$ and thus
$b_k \ge \Nm(\cI)$.
 Lemma~\ref{growth} now implies the desired bound.
\end{proof}

\section{Proof of Theorem~\ref{weakZ}}

For a prime ideal $\wp$ and an integer $k\ge 1$ we denote by
$\nu_{\wp}(k)$ the $\wp$-adic order of $k$. Let
$$
W = \prod_{n = M+1}^{M+N} \beta_n.
$$
For a prime ideal $\wp$, we denote by $r_{\wp}$ the
$\wp$-adic order of $W$
and by $s_{\wp}$ the largest $\wp$-adic order of the terms
$\beta_n$, $n = M+1, \ldots, M+N$.
Then
$$
r_{\wp} \le \sum_{s =1}^{s_{\wp}} T(M,N,\wp^s).
$$
By Lemma~\ref{propofT} we have
$$
r_{\wp} \ll \sum_{s =1}^{s_{\wp}} \left( \frac{N}{(s \log \Nm(\wp))^{1/2}}
+
1\right)
\ll  \frac{N s_{\wp}^{1/2}}{(\log \Nm(\wp))^{1/2}} + s_{\wp}.
$$
By Lemma~\ref{growth} we see that $s_{\wp} \ll (N+M)^2/\log \Nm(\wp)$
therefore
$$
r_{\wp} \ll  \frac{(N+M)^2}{\log \Nm(\wp)}.
$$
Thus
\begin{eqnarray*}
\log \Nm(W)  &= & \sum_{r_{\wp} > 0 } r_{\wp} \log \Nm(\wp) \ll
(N+M)^2
\sum_{r_{\wp} > 0 } 1 \\
  &= &  (N+M)^2  \omega_\K\left(\prod_{n = M+1}^{M+N}\beta_n\right).
\end{eqnarray*}
Lemma~\ref{growth} implies that $N (N+M)^2 \ll \log \Nm(W)$
which finishes the proof. \qed
\bigskip


\section{Proof of Theorem~\ref{better}}

\begin{lemma}\label{divseq}
If $Q= O$, then the sequence $\beta_n$ is a divisibility sequence, meaning
that $\beta_m|\beta_n$
as ideals, whenever $m|n$.
\end{lemma}
\begin{proof}
The proof of this follows
from the standard local theory of elliptic curves, see~\cite{SIL1}.
For every prime ideal $\wp$, write
$\K_{\wp}$ for the completion of $\K$ with respect to the
valuation corresponding to $\wp$. There is
a subgroup of the group of $\K_{\wp}$-rational
points:
$$\H(\K_{\wp})=\{O\} \cup \{R\in \E(\K_{\wp}):\mbox{ord} _{\wp}(x(R))
\leq -2\}.$$
Silverman~\cite{SIL1} gives a proof that for all $R\in \H(\K_{\wp})$,
$$\mbox{ord} _{\wp}(x(nR))=\mbox{ord} _{\wp}(x(R))-\mbox{ord} _{\wp}(n)$$
and
the divisibility statement follows at once from this.
\end{proof}
\begin{proof}[Proof of Theorem~\ref{better}.]
Suppose $r \geq  3w$ is an integer and $r < p_1 < \ldots < p_r$
denote $r$ distinct primes. Put $q=p_1\ldots p_r$.
Recalling the notation of~(\ref{definegenbeta}),
let $R=qP$ and write $(x(nR))=\gamma_n/\delta_n$
for the factorisation into integral ideals. Similarly,
for each $i=1,\ldots ,r$, write $R_i=p_iP$ with
$(x(nR_i))=\gamma_{in}/\delta_{in}$ for the
factorisation into integral ideals.
We claim
that if  $n$ is coprime to $q$, then $\delta_n$  has
at least $r$ distinct primitive prime ideal divisors.

From Cheon and Hahn's Theorem~\cite{ch}, for all sufficiently
large $n$, $\delta_{in}$ has a primitive prime ideal divisor. We
claim each of these is a primitive prime ideal divisor of
$\delta_n$ when $n$ is coprime to $q$. By Lemma~\ref{divseq},
$\beta$ is a divisibility sequence, thus any divisor of
$\delta_{in}$ is a divisor of $\delta_n$. Suppose $\wp$ is a
primitive prime ideal divisor of $\delta_{in}$ and $\wp$ divides
$\delta_m$ for some $m$. Then $mR=mqP\equiv O$ mod $\wp$. In other
words, $(mq/p_i)(p_iP)\equiv O$ mod $\wp$. Now $\wp$ is a
primitive divisor of $\delta_{in}$ and this forces $n$ to divide
$mq/p_i$. Since $n$ is coprime to $q$, this forces $n$ to divide
$m$ and shows that every primitive ideal divisor of $\delta_{in}$
is a primitive ideal divisor of $\delta_n$. Each of these
primitive ideal divisors is distinct thus we conclude that
$\delta_n$ has at least $r$ distinct primitive prime ideal
divisors. It follows that for $M$ sufficiently large,
\begin{eqnarray*}
\omega_\K \left( \prod_{n=M+1}^{M+N}\beta_{qn}\right)
& \geq& \omega_K \left( \prod_{n=M+1,(n,q)=1}^{M+N}\delta_n\right)
\geq rN\prod_{i=1}^r\left(1-\frac{1}{p_i}\right)\\
& \geq & rN \left(1-\frac{1}{r}\right)^r \ge re^{-1} N \ge wN,
\end{eqnarray*}
which  completes the proof of Theorem~\ref{better}.
\end{proof}

\section{Somos Sequences}

Sequences such as $u_n=a^n-b^n$ are known to satisfy a binary
linear recurrence relation with integer coefficients. Recently,
Swart~\cite{swart}, building on some unpublished identities of
Nelson Stephens, has related elliptic curves to {\it Somos sequences}.
These sequences, named after Michael Somos \cite{So},
satisfy a bilinear recurrence relation of the form
\begin{equation}\label{brr}
s_{n+2}s_{n-2}=As_{n+1}s_{n-1}+Bs_n^2,
\end{equation}
where $A$ and $B$ are rational numbers, not both zero. A special
case is the sequence $s$ with $A=B=1$, which begins $1,1,1,1,\ldots$ with
subsequent terms defined by the recurrence relation (\ref{brr}).
This is known as the {\it Somos-4 sequence}.
Although a division is needed in the definition of the later
terms, they do all turn out to be integral. For an introduction to
Somos sequences, as well as links to other areas of
mathematics, consult~\cite{epsw},
~\cite{dg1},~\cite{dg2},~\cite{jm},~\cite{Pr},
~\cite{rr}.

If $P=(x,y)$ and $Q$ are rational points with $Q+nP$
never the point at infinity, continue to write $Q+nP=(x_n,y_n)$.
Then a sequence $s$ satisfying (\ref{brr}) can be
obtained as follows (see~\cite{swart}): let $s_{-1}$
and $s_0$ be arbitrary non-zero rational numbers and
\begin{equation}\label{sformulae}
s_{n+1} = -\frac{(x_n - x)s_n^2}{s_{n-1}} \mbox{ for } n \geq 0.
\end{equation}

Clearly the sequence $s_n$ could turn out to have rational
terms. In fact Swart~\cite{swart} shows that the denominators are
easy to understand in terms of the starting data. In a
sense that can be made precise, each rational rational
sequence is equivalent to an integral sequence. Moreover,
she shows that beginning with a
rational sequence $s$
satisfying (\ref{brr}) it is possible to
reverse the process and
recover an elliptic curve together with the two rational points
$P$ and $Q$ which yields the sequence $s$ according to the
formula in (\ref{sformulae}).

\begin{example}
If $s$ is the Somos-4 sequence, the corresponding elliptic
curve is
$$y^2+y=x^3-x,
$$
with $Q=(0,0)$ and $P=(1,0)=2Q$ and $s_{-1}=s_0=1$.
Thus the terms of the Somos-4 sequence
correspond to the denominators of the odd multiples of
the point $Q=(0,0)$.
\end{example}

Any {\it elliptic divisibility sequence\/} (see~\cite{S,W1,W2} for
background and basic properties) arises when $Q$ is the point at
infinity, and satisfies a relation~(\ref{brr}) with $A=1$ and
$B=-1$. Thus Somos sequences are a generalisation of elliptic
divisibility sequences. We are able to record a weaker version of
Theorem \ref{weakZ} for Somos sequences. This is stated in the
rational case because it relies on Swart's thesis, which deals
with rational sequences only (although her methods will surely
generalise to algebraic sequences). The definition of $\omega$
extends to all non-zero rational numbers, if $q=a/b$ with $a$ and
$b$ denoting coprime integers then $\omega(q)=\omega(a)$.

\begin{corollary}\label{cor2}
Let $\E$ denote a rational elliptic curve having rational points
$P$ and $Q$ with $P$ non-torsion and $Q+nP=(x_n,y_n)$ not
the point at infinity for all $n$. Let $s$ denote an integral
sequence defined by (\ref{sformulae}), which
satisfies (\ref{brr}) for some $A$ and $B$. Then
$$
\omega\left(\prod_{n=1}^{N}s_n\right) \gg N.
$$
\end{corollary}

\begin{proof}
It is
shown in~\cite{ek} that for any Somos sequence $s$,
any primitive divisor of $b_n$
is a primitive divisor of the sequence
of denominators of $s_{n-1}$.
The result now follows from  Corollary~\ref{cor1}
\end{proof}

It seems likely that the methods
in this paper could be used to prove
a stronger version of Corollary~\ref{cor2} true
on intervals.


\begin{thebibliography}{9999}


\bibitem{bhv}
Yu. Bilu, G. Hanrot and P. M. Voutier,
`Existence of primitive divisors of Lucas and
Lehmer numbers'
{\it J. Reine Angew. Math.\/}, 539 (2001), 75--122.

\bibitem{ch}
J. Cheon and S. Hahn,
`The orders of the reductions of a point in the Mordell-Weil
group of an elliptic curve',
{\it Acta. Arith.\/}, {\bf 88} no. 3 (1999), 219--222.

\bibitem{dav}
S.  David,
`Minorations de formes lin\'eaires de logarithmes
elliptiques',
{\it M\'em. Soc. Math. France (N.S.)\/}, {\bf 62}, 1995.

\bibitem{eg}
G. Everest and K.  Gy\" ory,
`Primitive prime divisors', {\it Preprint\/}, 2003.

\bibitem{ek}
G. Everest and H. King,
`Primitive divisors of bilinear recurrence sequences',
{\it Preprint\/}, 2003.

\bibitem{epsw}
G. Everest, A. J. van der Poorten, I. E. Shparlinski
and T. Ward,
{\it Recurrence sequences\/},
Math. Surveys and Monographs, {\bf 104},
Amer. Math. Soc., Providence, RI, 2003.

\bibitem{dg1}
D. Gale,
`The strange and surprising saga of the Somos sequences',
{\it Math. Intelligencer\/}, {\bf 13} (1991),   40--42.

\bibitem{dg2}
D. Gale,
`Somos sequence update',
{\it Math. Intelligencer\/}, {\bf 13} 1991, 49--50.

\bibitem{L}
D. H. Lehmer,
`Factorization of certain cyclotomic functions',
{\it Ann. of Math.\/}, {\bf 34} (1933), 461--479.

\bibitem{jm}
J. L. Malouf,
`An integer sequence from a rational recursion',
{\it Discrete Math.\/}, {\bf 110} (1992),  257--261.

\bibitem{phi}
T. Pheidas,
`An effort to prove that the existential theory
of $\Q$ is undecidable',
{\it Hilbert's tenth problem: relations with arithmetic and
algebraic geometry (Ghent 1999)\/},
  Contemp. Math., {\bf 207}, Amer. Math. Soc., Providence, RI,
2000, 237--252.

\bibitem{P}
T. A. Pierce,
`Numerical factors of the arithmetical forms
$\prod_{i=1}^n(1\pm\alpha_i^m)$',
{\it Ann. of Math.\/}, {\bf 18} (1917), 53--64.

\bibitem{Pr}
J.  Propp,
`The Somos sequence site', Available from
{\tt http://www.math.wisc.edu/\~{ }propp/somos.html}.

\bibitem{rr}
R. M. Robinson,
`Periodicity of Somos sequences',
{\it Proc. Amer. Math. Soc.\/}, {\bf 116} (1992),  613--619.

\bibitem{sch}
A. Schinzel,
`On primitive prime factors of $a^n-b^n$',
{\it Proc. Camb. Phil. Soc.\/}, {\bf 58} (1962), 555--562.

\bibitem{Sch}
A. Schinzel,
`Primitive divisors of the expression $A^n-B^n$ in
algebraic number fields',
{\it J. reine angew. Math.\/}, {\bf 268/69} (1974), 27--33.

\bibitem{Sch2}
A. Schinzel,
`An extension of the theorem on primitive
divisors in
algebraic number fields',
{\it Math. Comp.\/}, {\bf 61} no. 203 (1993), 441--444.

\bibitem{S}
R. Shipsey,
`Elliptic divisibility sequences',
{\it PhD Thesis\/}, Univ. of London, 2000.

\bibitem{SIL1}
J. H. Silverman, {\it The {A}rithmetic of elliptic curves\/},
Springer, New York, 1986.

\bibitem{SIL3}
J. H. Silverman,
`Weiferich's criterion and the $ABC$-conjecture',
{\it J. Number Theory\/},  {\bf 30} (1988), 226--237.

\bibitem{SIL2}
J. H. Silverman, {\it Advanced topics in the arithmetic of
elliptic curves\/}, Springer, New York, 1994.


\bibitem{So}
M. Somos,
`Problem 1470',
{\it Crux Mathematicorum.\/},  {\bf 15} (1989), 208.

\bibitem{ST}
R. J. Stroeker and N. Tzanakis,
`Solving elliptic Diophantine equations by estimating
linear forms in elliptic logarithms',
{\it Acta. Arith.\/},  {\bf 67} (1994), 177--196.

\bibitem{swart}
C. Swart,
`Elliptic divisibility sequences', {\it PhD Thesis\/},
Univ. of London, 2003.

\bibitem{W1} M. Ward, `The law of repetition of primes in an elliptic
divisibility sequence',  {\it Duke Math. J.\/},  {\bf 15} (1948),
941--946.

\bibitem{W2} M. Ward, `Memoir on elliptic divisibility sequences',
{\it Amer. J. Math.\/},  {\bf 70} (1948) 31--74.

\bibitem{Z} K. Zsigmondy, `Zur Theorie der Potenzreste',
{\it Monatsh. Math.\/}, {\bf 3} (1892),  265--284.


\end{thebibliography}
\end{document}